\documentclass[12pt]{article}
\usepackage{color}
\definecolor{darkblue}{rgb}{0.00,0.25,0.50}
\usepackage[colorlinks,filecolor=blue,citecolor=darkblue]{hyperref}

\usepackage{amsfonts,amssymb,amsmath,amsthm}
\usepackage{url}
\usepackage{enumerate}
\usepackage[ukrainian, russian, english]{babel}
\usepackage[cp1251]{inputenc}
\begin{document}\selectlanguage{ukrainian}
\thispagestyle{empty}

\title{}

\begin{center}
\textbf{\Large Порядкові оцінки найкращих наближень і наближень
сумами Фур'є класів нескінченно диференційовних функцій}
\end{center}
\vskip0.5cm
\begin{center}
А.~С.~Сердюк${}^1$, Т.~А.~Степанюк${}^2$\\ \emph{\small
${}^1$Інститут математики НАН
України, Київ\\
${}^2$Східноєвропейський національний університет імені Лесі
Українки, Луцьк\\}
\end{center}
\vskip0.5cm


\begin{abstract}
Отримано порядкові оцінки для найкращих рівномірних наближень
тригонометричними поліномами та наближень сумами Фур'є класів
$2\pi$--періодичних неперервних функцій, таких, що їх
$(\psi,\beta)$--похідні $f_{\beta}^{\psi}$ належать одиничним кулям
просторів $L_{p}, \ 1\leq p<\infty$, у випадку коли послідовності
$\psi(k)$ спадають до нуля швидше за будь--яку степеневу функцію.
Аналогічні оцінки одержані для наближень в $L_{s}$--метриці, $1<
s\leq \infty$, для класів сумовних $(\psi,\beta)$--диференційовних
функцій, таких, що $\parallel f_{\beta}^{\psi}\parallel_{1}\leq1$.

\vskip 0.5cm

We obtained order estimations for the best uniform approximation by
trigonometric polynomials  and approximation by Fourier sums of
classes of  $2\pi$--periodic continuous functions, whose
$(\psi,\beta)$--derivatives $f_{\beta}^{\psi}$ belong to unit balls
of spaces $L_{p}, \ 1\leq p<\infty$ in case at consequences
$\psi(k)$ decrease to nought faster than any power function. We also
established the analogical estimations in $L_{s}$--metric, $1< s\leq
\infty$, for classes
 of the summable
$(\psi,\beta)$--differentiable functions, such that  $\parallel
f_{\beta}^{\psi}\parallel_{1}\leq1$.

\end{abstract}


Нехай $C$ --- простір $2\pi$--періодичних неперервних функцій, у
якому норма задана за допомогою рівності $
{\|f\|_{C}=\max\limits_{t}|f(t)|}$; $L_{\infty}$ --- простір
$2\pi$--періодичних вимірних і суттєво обмежених функцій $f(t)$ з
нормою $\|f\|_{\infty}=\mathop{\rm{ess}\sup}\limits_{t}|f(t)|$;
$L_{p}$, $1\leq p<\infty$, --- простір $2\pi$--періодичних сумовних
в $p$--му степені на $[0,2\pi)$ функцій $f(t)$, в якому  норма
задана формулою
$\|f\|_{p}=\Big(\int\limits_{0}^{2\pi}|f(t)|^{p}dt\Big)^{\frac{1}{p}}.$

Нехай далі функція $f(x)\in L_{1}$, і її ряд Фур'є має вигляд

$$
f(x)\sim\frac{a_{0}(f)}{2}+\sum\limits_{k=1}^{\infty}a_{k}(f)\cos
kx+b_{k}(f)\sin kx,
$$
$\psi(k)$ --- фіксована послідовність дійсних чисел, а $\beta$
--- деяке  дійсне число. Якщо ряд
$$
\sum_{k=1}^{\infty}\frac{1}{\psi(k)}\left(a_{k}(f)\cos\left(kx+\frac{\beta\pi}{2}\right)
+b_{k}(f)\sin\left(kx+\frac{\beta\pi}{2}\right)\right)
$$
\noindent є рядом Фур'є деякої сумовної функції $\varphi$, то
функцію $\varphi$ називають (див., наприклад, \cite[с.
132]{Stepanets1}) $(\psi,\beta)$-похідною функції $f(x)$ і
позначають $f_{\beta}^{\psi}(x)$.

\rm Множину функцій $f(x)$, у яких існує $(\psi,\beta)$-похідна
позначають через $L_{\beta}^{\psi}$, а підмножину неперервних
функцій із $L_{\beta}^{\psi}$ -- через $C_{\beta}^{\psi}$.

Якщо $f\in L^{\psi}_{\beta}$, і водночас $f^{\psi}_{\beta}\in
\mathfrak{N}\subseteq L_{1}$, то кажуть, що функція $f$ належить
класу $L^{\psi}_{\beta}\mathfrak{N}$. Якщо $\mathfrak{N}=B_{p}^{0}$,
де ${\mathfrak{N}=B_{p}^{0}=\left\{\varphi: \ ||\varphi||_{p}\leq 1,
\ \varphi\perp1\right\}}$, то покладають
  $L^{\psi}_{\beta}B_{p}^{0}=L^{\psi}_{\beta,p}$.

  Як показано в \cite[с.
136]{Stepanets1}, якщо послідовність $\psi(k)$
  монотонно прямує до нуля при  $k\rightarrow\infty$
  i $\sum\limits_{k=1}^{\infty}\frac{\psi(k)}{k}<\infty$, то елементи $f(\cdot)$ множини
  $L^{\psi}_{\beta}\mathfrak{N}$ при будь--якому $\beta\in\mathbb{R}$ і майже для всіх $x\in\mathbb{R}$ представляються згортками
 \begin{equation}\label{psi,beta}
f(x)=\frac{a_{0}(f)}{2}+\frac{1}{\pi}\int\limits_{-\pi}^{\pi}\Psi_{\beta}(x-t)\varphi(t)dt,
\ a_{0}\in\mathbb{R}, \  \varphi\in\mathfrak{N}, \ \varphi\perp 1,
\end{equation}
з сумовним ядром $\Psi_{\beta}(t)$  ряд Фур'є якого має вигляд
 $$
\sum\limits_{k=1}^{\infty}\psi(k)\cos
\big(kt-\frac{\beta\pi}{2}\big).
$$
При цьому функція $\varphi(\cdot)$ майже скрізь співпадає з
$f^{\psi}_{\beta}(\cdot)$. Якщо ж $f\in
C^{\psi}_{\beta}\mathfrak{N}, \ \beta\in\mathbb{R}$, то рівність
(\ref{psi,beta}) виконується для всіх $x\in\mathbb{R}$.

Якщо $\psi(k)=k^{-r},\ r>0$, то класи $L_{\beta,p}^{\psi}$ є
відомими класами
 Вейля-Надя $W_{\beta,p}^r$.

Будемо вважати, що послідовності $\psi(k)$ є звуженням на множину
натуральних чисел $\mathbb{N}$ деяких додатних, неперервних, опуклих
донизу функцій $\psi(t)$ неперервного аргументу $t\geq1$ таких, що $
\lim\limits_{t\rightarrow\infty}\psi(t)=0. $ Множину всіх таких
функцій $\psi(t)$ позначатимемо через ${\mathfrak M}$.

Згідно з \cite[с. 159--160]{Stepanets1}, кожній функції
$\psi\in{\mathfrak M}$ поставимо у відповідність характеристики
$$
\eta(t)=\eta(\psi;t)=\psi^{-1}\left(\psi(t)/2\right), \ \ \
\mu(t)=\mu(\psi;t)=\frac{t}{\eta(t)-t},
$$
де $\psi^{-1}(\cdot)$ --- обернена до $\psi$ функція і розглянемо
множину
$$
\mathfrak{M}^{+}_{\infty}=\left\{\psi\in \mathfrak{M}: \ \
\mu(\psi;t)\uparrow\infty, \ t\rightarrow\infty \right\}.
$$

Через ${\mathfrak M^{'}_{\infty}}$ позначимо підмножину функцій
$\psi\in {\mathfrak M^{+}_{\infty}}$, для кожної з яких величина
$\eta(\psi;t)-t$ обмежена зверху, тобто існує стала ${K_{1}>0}$,
така, що $\eta(\psi;t)-t\leq K_{1}, \ t\geq1$, а через ${\mathfrak
M^{''}_{\infty}}$ --- підмножину функцій ${\psi\in {\mathfrak
M^{+}_{\infty}}}$, для кожної з яких величина $\eta(\psi;t)-t$
обмежена знизу деяким додатним числом, тобто існує стала $K_{2}>0$,
така, що ${\eta(\psi;t)-t\geq K_{2}}, \ t\geq1$.

Типовими представниками множини $\mathfrak{M}^{+}_{\infty}$  є
функції ${\psi_{r}(t)=\exp(-\alpha t^{r})}, \ {\alpha>0}, \ r>0$,
причому, якщо $r\geq1$, то $\psi_{r}\in \mathfrak{M}^{'}_{\infty}$,
а якщо ${r\in(0,1]}$, то $\psi_{r}\in \mathfrak{M}^{''}_{\infty}$.

Через $F$ прийнято позначати множину функцій $\psi\in {\mathfrak
M}$, таких що $\eta'(\psi,t)=\eta'(\psi,t+0)\leq K$. Відмітимо
(див., наприклад, \cite[с. 165]{Stepanets1}), що
$\mathfrak{M}^{+}_{\infty}\subset F$.

Якщо $\psi\in\mathfrak{M}^{+}_{\infty}$, то (див., наприклад,
\cite[с. 97]{Step monog 1987}) множини $C^{\psi}_{\beta}$
складаються з нескінченно диференційовних функцій. З іншого боку, як
показано в \cite[с. 1692]{Stepanets_Serdyuk_Shydlich}, для кожної
нескінченно диференційовної, $2\pi$--періодичної функції $f$ можна
вказати функцію $\psi$ з множини $\mathfrak{M}^{+}_{\infty}$, таку
що $f\in C^{\psi}_{\beta}$ для довільних $\beta\in\mathbb{R}$.

Метою даної роботи є  знаходження точних порядкових оцінок для
 величин вигляду
$$
{\cal E}_{n}(\mathfrak{N})_{X}=\sup\limits_{f\in
\mathfrak{N}}\|f(\cdot)-S_{n-1}(f;\cdot)\|_{X},
$$
де $S_{n-1}(f;\cdot)$ --- частинні суми Фур'є порядку $n-1$,
$\mathfrak{N}\subset X\subset L_{1}$, а також знаходження точних
порядкових оцінок найкращих наближень, тобто величин вигляду
$$
{E}_{n}(\mathfrak{N})_{X}=\sup\limits_{f\in
\mathfrak{N}}\inf\limits_{t_{n-1}\in\mathcal{T}_{2n-1}}\|f(\cdot)-t_{n-1}(\cdot)\|_{X},
$$
де $\mathcal{T}_{2n-1}$ --- підпростір усіх тригонометричних
поліномів $t_{n-1}$ порядку не вищого за $n-1$, у наступних
випадках:

1) $\mathfrak{N}=C^{\psi}_{\beta,p}, \ 1\leq p<\infty, \ X=C$;

2) $\mathfrak{N}=L^{\psi}_{\beta,1}, \ X=L_{s}, \ 1<s\leq\infty$

\noindent при $\psi\in{\mathfrak M^{''}_{\infty}}$ і
$\beta\in\mathbb{R}$.

При $X=L_{s}, 1\leq s\leq \infty$ будемо позначати ${\cal
E}_{n}(\mathfrak{N})_{L_{s}}$ через ${\cal
E}_{n}(\mathfrak{N})_{s}$, а $E_{n}(\mathfrak{N})_{L_{s}}$ через
${\cal E}_{n}(\mathfrak{N})_{s}$ відповідно.

Зробимо короткий історичний огляд  дослідження величин ${\cal
E}_{n}(L^{\psi}_{\beta,p})_{s}$ і $E_{n}(L^{\psi}_{\beta,p})_{s}$.

Для класів Вейля--Надя   $W^{r}_{\beta,p}$, при довільних $r>0$,
$\beta\in\mathbb{R}$, ${1\leq p,s\leq\infty}$ точні порядкові оцінки
величин ${\cal E}_{n}(W^{r}_{\beta,p})_{{s}}$,
${E}_{n}(W^{r}_{\beta,p})_{{s}}$ відомі (див., наприклад, \cite[с.
47--49]{T}).

У випадку $p=s=1$ i $p=s=\infty$ відомі також асимптотичні рівності
при $n\rightarrow\infty$ для величин ${\cal
E}_{n}(W^{r}_{\beta,\infty})_{\infty}$ та ${\cal
E}_{n}(W^{r}_{\beta,1})_{1} $, $r>0$, $\beta\in\mathbb{R}$ (див.,
наприклад,  роботи \cite{Kol}--\cite{Nik2}).

У випадках $p=s=1$,  $p=s=\infty$, $p=s=2$ та $p=2$ і $s=\infty$
встановлені точні значення  найкращих наближень
${E}_{n}(W^{r}_{\beta,\infty})_{\infty}$ та
${E}_{n}(W^{r}_{\beta,1})_{1}$  при усіх $n\in\mathbb{N}$, $r>0$ і
$\beta\in \mathbb{R}$ (див.  роботи \cite{Fav}--\cite{Bab}.

На класах $L^{\psi}_{\beta,p}$ точні порядкові оцінки величин ${\cal
E}_{n}(L^{\psi}_{\beta,p})_{s}$ та $E_{n}(L^{\psi}_{\beta,p})_{s}$ у
випадку, коли $\psi(k)k^{\frac{1}{p}-\frac{1}{s}}$ монотонно
незростають і $\psi\in B$, де  $B$ --- множина монотонно
незростаючих при $t\geq 1$ додатних функцій $\psi(t)$, для кожної з
яких можна вказати додатну сталу $K$ таку, що $
\frac{\psi(t)}{\psi(2t)}\leq K, \ \  t\geq 1 $, були знайдені у
роботі \cite{Kuchpel}  при довільних ${1< p,s<\infty}$.

При $p=s=2$ в \cite{Kuchpel} також розв'язано задачу про точні
значення величин ${\cal E}_{n}(L^{\psi}_{\beta,p})_{s}$ та
$E_{n}(L^{\psi}_{\beta,p})_{s}$ за умови $\sup\limits_{k\geq
n}\psi(k)<\infty$.

 Зазначимо також, що при $p=2$ і $s=\infty$ або $p=1$ і $s=2$ точні значення величин
${\cal E}_{n}(C^{\psi}_{\beta,p})_{s}$  для всіх $n\in\mathbb{N}$,
$\beta\in \mathbb{R}$
 за умови збіжності ряду
$\sum\limits_{k=1}^{\infty}\psi^{2}(k)$ знайдені у роботах
\cite{Serdyuk} і \cite{Serdyuk2013}.

В \cite{Serdyuk_grabova}  встановлено  точні порядкові оцінки
величин ${\cal E}_{n}(L^{\psi}_{\beta,p})_{s}$ та
$E_{n}(L^{\psi}_{\beta,p})_{s}$ при $1\leq p<\infty, \ s=\infty$, а
також при $p=1$ i $1<s\leq\infty$, у випадку коли $\psi\in
B\cap\Theta_{p}$, де $\Theta_{p}$, $1\leq p<\infty$, --- множина
монотонно незростаючих функцій $\psi(t)$, для яких існує стала
$\alpha>\frac{1}{p}$ така, що функція $t^{\alpha}\psi(t)$ майже
спадає.

При $\psi\in {\mathfrak M^{'}_{\infty}}$, для довільних $1\leq
p,s\leq\infty$, $n\in\mathbb{N}$ і
 $\beta\in\mathbb{R}$ точні порядкові
оцінки величин ${\cal E}_{n}(L^{\psi}_{\beta,p})_{s}$ та
$E_{n}(L^{\psi}_{\beta,p})_{s}$, встановлені в  \cite[с. 225]{Step
monog 1987}    (див. також  \cite[с. 48]{Stepanets2}) і мають вигляд
\begin{equation}\label{rezm1}
C^{(1)}_{p,s}\psi(n)\leq E_{n}(L^{\psi}_{\beta,p})_{s}\leq {\cal
E}_{n}(L^{\psi}_{\beta,p})_{s}\leq C^{(2)}_{p,s}\psi(n),
\end{equation}
де $C^{(1)}_{p,s}, \ C^{(2)}_{p,s}$ --- додатні сталі, що залежать
тільки від $p$ i $s$.

В  \cite[с. 219]{Step monog 1987}   (див. також \cite[с.
60]{Stepanets2}) знайдено також точні порядкові оцінки величин
${\cal E}_{n}(L^{\psi}_{\beta,p})_{s}$ та
$E_{n}(L^{\psi}_{\beta,p})_{s}$,
 при $\psi\in {\mathfrak
M^{''}_{\infty}}$, для довільних $1<p,s<\infty$, $n\in\mathbb{N}$ і
$\beta\in\mathbb{R}$, які мають вигляд
\begin{equation}\label{rezSt}
C^{(3)}_{p,s}\psi(n)(\eta(n)-n)^{\alpha}\leq
E_{n}(L^{\psi}_{\beta,p})_{s}\leq {\cal
E}_{n}(L^{\psi}_{\beta,p})_{s}\leq
C^{(4)}_{p,s}\psi(n)(\eta(n)-n)^{\alpha},
\end{equation}
де $C^{(3)}_{p,s}, \ C^{(4)}_{p,s}$ --- додатні сталі, що залежать
тільки від $p$ i $s$, a ${\alpha=p^{-1}-s^{-1}}$, якщо $p<s$, i
$\alpha=0$, якщо $p\geq s$.

Оцінка зверху в співвідношенні (\ref{rezSt}) є справедливою і у
випадку $p=1$, за умови ${p<s<\infty}$ (див. \cite[с. 224]{Step
monog 1987}).

У випадках $p=s=1$,  $p=s=\infty$, $\psi\in {\mathfrak
M^{+}_{\infty}}$ і $\beta\in\mathbb{R}$ в \cite{serdyuk2004zbirnyk}
встановлено асимптотичні рівності для величин
$E_{n}(C^{\psi}_{\beta,p})_{s}$ і $E_{n}(L^{\psi}_{\beta,p})_{s}$.
Крім того в  \cite{Serdyuk2002}
  при $p=s=\infty$  та ${p=s=1}$ отримано  точні значення величин
$E_{n}(C^{\psi}_{\beta,p})_{s}$ i $E_{n}(L^{\psi}_{\beta,p})_{s}$,
${\beta\in\mathbb{R}}$ за умови, що функція $\psi(k),
{k\in\mathbb{N}}$ має наступні властивості: 1)
${\Delta^{2}\psi(k)\mathop{=}\limits^{\rm
df}\psi(k)-2\psi(k+1)+\psi(k+2)\geq 0}$,
${\frac{\psi(k+1)}{\psi(k)}\leq \rho}, \ {0<\rho<1}, \ {k=n,
n+1,...}$; 2)
${\frac{\Delta^{2}\psi(n)}{\psi(n)}>\frac{(1+3\rho)\rho^{2n}}{(1-\rho)\sqrt{1-2\rho^{2n}}}
} $.

В \cite{Rom} для ${\psi\in {\mathfrak M^{+}_{\infty}}}$, за умови,
що починаючи з деякого $t_{0}\geq 1$ ${\eta(t)-t> 1}$,
 встановлено
точні порядкові оцінки величин ${\cal
E}_{n}(C^{\psi}_{\beta,p})_{s}$,
 $\beta\in \mathbb{R}$ при ${1<p<\infty}, \ s=\infty$, котрі мають вигляд:
$$
C_{\psi,p}^{(1)} \psi(n)(\eta(n)-n)^{\frac{1}{p}}\leq{\cal
E}_{n}\Big(C^{\psi}_{\beta,p}\Big)_{C}\leq C_{\psi,p}^{(2)}
\psi(n)(\eta(n)-n)^{\frac{1}{p}},\ n\in\mathbb{N},
$$
де $C_{\psi,p}^{(1)}, \ C_{\psi,p}^{(2)}$ --- додатні сталі, що
залежать
 від $\psi$ i $p$.

В даній роботі встановлено точні порядкові оцінки величин
 ${E}_{n}\big(C^{\psi}_{\beta,p}\big)_{C}$,
${E}_{n}\big(L^{\psi}_{\beta,1}\big)_{s}$ i ${\cal
E}_{n}\big(L^{\psi}_{\beta,1}\big)_{s}$  для довільних ${1\leq
p<\infty}$, ${1<s\leq\infty}$ і $\beta\in\mathbb{R}$, у випадку,
коли $\psi\in\mathfrak{M}^{+}_{\infty}$,
$\lim\limits_{t\rightarrow\infty}(\eta(\psi,t)-t)=\infty$,
${\eta(t)-t\geq a>2}, \ {\mu(t)\geq b>2}$. При цьому константи в
порядкових оцінках записуються через параметри задачі в явному
вигляді.
 Отримані оцінки доповнюють і уточнюють згадані вище результати робіт \cite[с. 219]{Step
monog 1987}
 (див. також \cite[с.
60]{Stepanets2}) та  \cite{Rom}.

Перейдемо до викладу основних результатів.

\textbf{Теорема 1.} \emph{ Нехай $\psi\in\mathfrak{M}^{+}_{\infty},
\ \lim\limits_{t\rightarrow\infty}(\eta(\psi,t)-t)=\infty$,
$\beta\in \mathbb{R}$, ${1\leq p<\infty}$. Тоді  для  $n\in
 \mathbb{N}$, таких, що ${\eta(n)-n\geq a>2}, \ {\mu(n)\geq b>2}$ справедливі оцінки}
$$
C_{a}\psi(n)(\eta(n)-n)^{\frac{1}{p}} \leq {
E}_{n}\Big(C^{\psi}_{\beta,p}\Big)_{C}\leq{\cal
E}_{n}\Big(C^{\psi}_{\beta,p}\Big)_{C} \leq \ \ \ \ \ \ \ \ \ \ \ \
\ \ \ \ \ \ \ \ \ \ \ $$
\begin{equation}\label{theorem1}
\ \ \ \ \ \ \ \ \ \ \ \ \ \ \ \ \ \ \ \ \ \ \ \ \ \ \ \ \ \ \ \ \ \
\ \ \ \ \ \ \ \ \ \ \ \ \ \ \ \ \ \ \leq C_{a,b} \
(2p)^{1-\frac{1}{p}} \psi(n)(\eta(n)-n)^{\frac{1}{p}},
\end{equation}
де
\begin{equation}\label{Ca}
C_{a}=\frac{\pi}{96\left(1+\pi^{2}\right)^{2}}
\frac{(a-1)^{2}(a-2)^{2}}{a^{3}(3a-4)},
\end{equation}
\begin{equation}\label{Cab}
C_{a,b}= \frac{1}{\pi}\max\left\{\frac{2b}{b-2}+\frac{1}{a}, \
2\pi\right\}.
\end{equation}

\textbf{\emph{Доведення теореми.}}
 Спочатку оцінимо зверху величину ${\cal
 E}_{n}\Big(C^{\psi}_{\beta,p}\Big)_{C}$.
Згідно з інтегральним зображенням (\ref{psi,beta}), для довільної
функції $f\in L^{\psi}_{\beta,p}, \ 1\leq p\leq\infty$,
$\psi\in\mathfrak{M}^{+}_{\infty}$ майже скрізь виконується рівність
\begin{equation}\label{2t1}
f(x)-S_{n-1}(f;x)=\frac{1}{\pi}\int\limits_{-\pi}^{\pi}\Psi_{\beta,n}(x-t)\varphi(t)dt,
\end{equation}
 де
 \begin{equation}\label{t13}
\|\varphi\|_{p}\leq1, \ \varphi\perp1,
\end{equation}
\begin{equation}\label{2t2}
\Psi_{\beta,n}(t)=
\sum\limits_{k=n}^{\infty}\psi(k)\cos\big(kt-\frac{\beta\pi}{2}\big).
\end{equation}

При цьому, якщо $f\in C^{\psi}_{\beta,p}, \ 1\leq p\leq\infty$, то
рівність (\ref{2t1}) виконується в кожній точці.

Далі нам буде корисним твердження роботи \cite[с.
137--138]{Stepanets1}.

\textbf{Твердження 1.} \emph{Якщо $h\in L_{p}, \ 1\leq p\leq\infty$,
${g\in L_{p'}}, \ {\frac{1}{p}+\frac{1}{p'}=1}$, то згортка
 $$
 f(x)=\frac{1}{\pi}\int\limits_{\pi}^{\pi}h(x-t)g(t)dt
 $$
 неперервна на всій осі, причому}
\begin{equation}\label{statement1}
\left\|f\right\|_{C}\leq\frac{1}{\pi}\left\|h\right\|_{p}\left\|g\right\|_{p'}.
\end{equation}
В силу твердження 1, та формул (\ref{2t1}) i (\ref{t13})
 \begin{equation}\label{gg10}
{\cal E}_{n}\big(C^{\psi}_{\beta,p}\big)_{C}\leq
\frac{1}{\pi}\big\|\Psi_{\beta,n}(\cdot)\big\|_{p'}\|
\varphi(\cdot)\|_{p}\leq\frac{1}{\pi}\big\|\Psi_{\beta,n}(\cdot)\big\|_{p'},
\end{equation}
де $p'=\frac{p}{p-1}$.

Застосувавши до функції $\Psi_{\beta,n}(t)$ перетворення Абеля, при
довільному $n\in~ \mathbb{N}$ одержимо
\begin{equation}\label{3t2}
\Psi_{\beta,n
}(t)=\sum\limits_{k=n}^{\infty}(\psi(k)-\psi(k+1))D_{k,\beta}(t)-
\psi(n)D_{n-1,\beta}(t),
\end{equation}
де
$$
D_{k,\beta}(t)=\frac{1}{2}\cos\frac{\beta\pi}{2}+\sum\limits_{j=1}^{k}\cos\left(jt-\frac{\beta\pi}{2}\right).
$$
Врахування відомих формул (див., наприклад, \cite[с. 40,42]{Step
monog 1987})
\begin{equation}\label{Dir ker}
D_{k,0}(t)=\frac{1}{2}+\sum\limits_{j=1}^{k}\cos
kt=\frac{\sin\left(k+\frac{1}{2}\right)t}{2\sin\frac{t}{2}}, \
0<|t|\leq \pi,
\end{equation}
і
\begin{equation}\label{con Dir ker}
D_{k,1}(t)=\sum\limits_{j=1}^{k}\sin
kt=\frac{\cos\frac{t}{2}-\cos\left(k+\frac{1}{2}\right)t}{2\sin\frac{t}{2}},
\ 0<|t|\leq \pi,
\end{equation}
де $D_{k,0}$ --- ядро Діріхле порядку $k$, а $D_{k,1}$ --- спряжене
ядро Діріхле порядку $k$, дозволяє записати
$$
D_{k,\beta}(t)=\cos\frac{\beta\pi}{2}D_{k,0}(t)+\sin\frac{\beta\pi}{2}D_{k,1}(t)=
$$
$$
=\cos\frac{\beta\pi}{2}\frac{\sin\left(k+\frac{1}{2}\right)t}{2\sin\frac{t}{2}}+\sin\frac{\beta\pi}{2}
\frac{\cos\frac{t}{2}-\cos\left(k+\frac{1}{2}\right)t}{2\sin\frac{t}{2}}=
$$
\begin{equation}\label{gg8}
=\frac{\sin\left(\left(k+\frac{1}{2}\right)t-
\frac{\beta\pi}{2}\right)+\sin\frac{\beta\pi}{2}\cos\frac{t}{2}}{2\sin\frac{t}{2}},
\ 0<|t|\leq \pi.
\end{equation}
Оскільки
\begin{equation}\label{ineq}
\sin \frac{t}{2}\geq\frac{t}{\pi}, \ \
 0\leq t\leq\pi,
\end{equation}
то з (\ref{gg8}) одержимо
\begin{equation}\label{gg9}
|D_{k,\beta}(t)|\leq\frac{\pi}{|t|}, \ 0<|t|\leq \pi.
\end{equation}
З формул (\ref{3t2}) i (\ref{gg9}) випливає, що
\begin{equation}\label{xx}
\left|\Psi_{\beta,n}(t)\right|\leq 2\pi\psi(n)\frac{1}{|t|}, \
0<|t|\leq \pi.
\end{equation}

Крім того, згідно з (\ref{2t2}) для довільних $t\in\mathbb{R}$
\begin{equation}\label{t1}
\left|\Psi_{\beta,n}(t)\right|\leq
\sum\limits_{k=n}^{\infty}\psi(k)\leq
\psi(n)+\int\limits_{n}^{\infty}\psi(u)du.
\end{equation}

Для оцінки інтеграла в правій частині формули (\ref{t1}),
скористаємось наступним твердженням роботи \cite[с.
500]{Serdyuk2004}.

\textbf{Твердження 2.} \emph{Якщо функція
 $\psi\in\mathfrak{M}^{+}_{\infty}$, то для  довільного ${m\in\mathbb{N}}$,
 такого, що $\mu(\psi,m)>2$
 виконується умова}
\begin{equation}\label{statement2}
\ \ \int\limits_{m}^{\infty}\psi(u)du\leq
\frac{2}{1-\frac{2}{\mu(m)}}\psi(m)(\eta(m)-m).
\end{equation}

Якщо $\mu(\psi,n)\geq b>2$, то з  нерівності (\ref{statement2}),
отримуємо
\begin{equation}\label{st}
\ \ \int\limits_{n}^{\infty}\psi(u)du\leq
\frac{2b}{b-2}\psi(n)(\eta(n)-n).
\end{equation}
З урахуванням формул (\ref{t1}) i (\ref{st}), а також умови
$\eta(n)-n\geq a>0$, маємо
\begin{equation}\label{z10}
\left|\Psi_{\beta,n}(t)\right|\leq
\left(\frac{2b}{b-2}+\frac{1}{a}\right)\psi(n)(\eta(n)-n), \ a>0, \
b>2.
\end{equation}

Поклавши ${C_{a,b}=\frac{1}{\pi}\max\{\frac{2b}{b-2}+\frac{1}{a}, \
2\pi\}}$, та використовуючи формули (\ref{xx}), (\ref{z10}),
отримуємо
$$
\frac{1}{\pi}\|\Psi_{\beta,n}(t)\|_{p'}\leq
$$
$$
\leq C_{a,b}\psi(n)\left(\int\limits_{|t|\leq\frac{1}{\eta(n)-n}
}(\eta(n)-n)^{p'}dt+\int\limits_{\frac{1}{\eta(n)-n}\leq|t|\leq\pi}\frac{dt}{|t|^{p'}}\right)^{\frac{1}{p'}}\leq
$$
$$
\leq
C_{a,b}\psi(n)\left(\eta(n)-n\right)^{\frac{1}{p}}2^{\frac{1}{p'}}\left(1\!+\!
\frac{1}{p'-1}\left(1\!-\!\frac{\pi^{1-p'}}{(\eta(n)-n)^{p'-1}}\right)\!\!\right)^{\frac{1}{p'}}\!\!<
$$
$$
<C_{a,b} 2^{\frac{1}{p'}}\left(1+
\frac{1}{p'-1}\right)^{\frac{1}{p'}}\psi(n)\left(\eta(n)-n\right)^{\frac{1}{p}}=
$$
\begin{equation}\label{yy}
=
C_{a,b}2^{\frac{1}{p'}}p^{\frac{1}{p'}}\psi(n)\left(\eta(n)-n\right)^{\frac{1}{p}},
\ \ 1< p'<\infty, \ a>0, \ b>2.
\end{equation}

З рівності (\ref{z10}), отримуємо, що при $p'=\infty$
$$
\frac{1}{\pi}\big\|\Psi_{\beta,n}(\cdot)\big\|_{p'}=\frac{1}{\pi}\big\|\Psi_{\beta,n}(\cdot)\big\|_{\infty}\leq
\frac{1}{\pi}\left(\frac{2b}{b-2}+\frac{1}{a}\right)\psi(n)(\eta(n)-n)\leq
$$
\begin{equation}\label{q1}
\leq C_{a,b}\psi(n)(\eta(n)-n), \ a>0, \ b>2.
\end{equation}

Зі співвідношень   (\ref{gg10}), (\ref{yy}) i (\ref{q1}) при $a>0, \
b>2$, $1\leq p<\infty$ отримуємо оцінку
$$
{\cal E}_{n}\Big(C^{\psi}_{\beta,p}\Big)_{C}\leq
C_{a,b}(2p)^{1-\frac{1}{p}} \psi(n)(\eta(n)-n)^{\frac{1}{p}}.
$$

Для завершення доведення теореми 1, враховуючи очевидну нерівність
$$
{ E}_{n}\Big(C^{\psi}_{\beta,p}\Big)_{C}\leq{\cal
E}_{n}\Big(C^{\psi}_{\beta,p}\Big)_{C},
$$
потрібно показати, що за умови виконання умов
$\psi\in\mathfrak{M}^{+}_{\infty}$ i $\eta(n)-n\geq \linebreak\geq
a>2, \ \mu(n)\geq b>2$, знайдеться функція $f^{*}\in
C^{\psi}_{\beta,p}$, така, що
$$
{
E}_{n}(f^{*})_{C}=\inf\limits_{t_{n-1}\in\mathcal{T}_{2n-1}}\|f^{*}(\cdot)-t_{n-1}(\cdot)\|_{C}\geq
$$
\begin{equation}\label{t10}
\geq \frac{\pi}{96\left(1+\pi^{2}\right)^{2}}
\frac{(a-1)^{2}(a-2)^{2}}{a^{3}(3a-4)}\psi(n)(\eta(n)-n)^{\frac{1}{p}},
\ 1\leq p\leq\infty.
\end{equation}

Позначивши цілу частину дійсного числа $\alpha$ через $[\alpha]$,
розглянемо при заданому $n\in \mathbb{N}$ функцію
$$
f_{p}(t)=f_{p}(\psi;n;t)=
\frac{(a-1)(a-2)}{2\left(1+\pi^{2}\right)a(3a-4)}\frac{1}{(\eta(n)-n)^{1-\frac{1}{p}}}\times
$$
\begin{equation}\label{function}
\times\left(W_{[\eta(n)],[\eta(\eta(n))]}(\psi;0;t)-
W_{n,[\eta(n)]}(\psi;0;t)\right), \ a>2,
\end{equation}
в якій величини $W_{N,M}(\lambda;\gamma;t)$, $N,M\in\mathbb{N} \
(N<M)$ означаються за допомогою формули
\begin{equation}\label{www}
W_{N,M}(\lambda;\gamma;t)=\frac{1}{M-N}\sum\limits_{k=N}^{M-1}
\sum\limits_{j=1}^{k}\lambda(j)\cos\left(jt+\gamma\right),
\end{equation}
де $\gamma\in\mathbb{R}$, а  ${\lambda(k)}, \ {k=1,2,...}$ ---
фіксована послідовність дійсних чисел.

Покажемо спочатку, що функція $f_{p}(\cdot)$ належить класу
$C^{\psi}_{\beta,p}, {1\leq p\leq \infty}$. Для цього досить
переконатись у виконанні нерівності
\begin{equation}\label{gg11}
\Big\|\big(f_{p}(\cdot)\big)^{\psi}_{\beta}\Big\|_{p}\leq 1, \ 1\leq
p\leq  \infty.
\end{equation}
З цією метою спочатку покажемо, що для функції
$W_{N,M}(\lambda;\gamma;t)$ справедливе твердження.

\textbf{Лема 1.} \emph{Нехай $\gamma\in\mathbb{R}$, а ${\lambda(k)},
\ {k=1,2,...}$ --- деяка послідовність дійсних чисел. Тоді для
довільних $N,M\in\mathbb{N} \ (N<M)$ має місце рівність}
$$
W_{N,M}(\lambda;\gamma;t)=\sum\limits_{k=1}^{N}\lambda(k)\cos\left(kt+\gamma\right)+
$$
\begin{equation}\label{lemma1}
+
\frac{1}{M-N}\sum\limits_{k=N+1}^{M-1}\left(M-k\right)\lambda(k)\cos\left(kt+\gamma\right).
\end{equation}
 \emph{Доведення леми 1.} Рівність (\ref{lemma1}) випливає із наступного ланцюжка перетворень:
 $$
 \frac{1}{M-N}\sum\limits_{k=N}^{M-1}\sum\limits_{j=1}^{k}\lambda(j)\cos\left(jt+\gamma\right)
=
 $$
 $$
 =\frac{1}{M-N}\left(\sum\limits_{j=1}^{N}\lambda(j)\cos\left(jt+\gamma\right)+...+
 \sum\limits_{j=1}^{M-1}\lambda(j)\cos\left(jt+\gamma\right)\right)=
 $$
 $$
 =\sum\limits_{k=1}^{N}\lambda(k)\cos\left(kt+\gamma\right)+
 $$
 $$
 +\frac{1}{M-N}\left[(M-N-1)\lambda(N+1)\cos\left((N+1)t+\gamma\right)+...+
 \right.
 $$
 $$
 \left.+\lambda(M-1)\cos\left((M-1)t+\gamma\right)\right]=\sum\limits_{k=1}^{N}\lambda(k)\cos\left(kt+\gamma\right)+
 $$
 $$
 +
\frac{1}{M-N}\sum\limits_{k=N+1}^{M-1}\left(M-k\right)\lambda(k)\cos\left(kt+\gamma\right).
 $$
 Лему 1 доведено.

Двічі використавши рівність (\ref{lemma1}), отримуємо
$$
W_{[\eta(n)],[\eta(\eta(n))]}(\psi;0;t)-
W_{n,[\eta(n)]}(\psi;0;t)=\sum\limits_{k=1}^{[\eta(n)]}\psi(k)\cos
kt+
$$
$$
+\frac{1}{[\eta(\eta(n))]-[\eta(n)]}
\sum\limits_{k=[\eta(n)]+1}^{[\eta(\eta(n))]-1}\left([\eta(\eta(n))]-k\right)\psi(k)\cos
kt-
$$
$$
-\sum\limits_{k=1}^{n}\psi(k)\cos kt-
\frac{1}{[\eta(n)]-n}\sum\limits_{k=n+1}^{[\eta(n)]-1}
\left([\eta(n)]-k\right)\psi(k)\cos kt =
$$
$$
=\sum\limits_{k=n+1}^{[\eta(n)]}\psi(k)\cos kt-
\frac{1}{[\eta(n)]\!-\!n}\!\!\sum\limits_{k=n+1}^{[\eta(n)]-1}\!\!\!
\left([\eta(n)]-k\right)\psi(k)\cos kt+
$$
$$
+\frac{1}{[\eta(\eta(n))]-[\eta(n)]}
\sum\limits_{k=[\eta(n)]+1}^{[\eta(\eta(n))]-1}\left([\eta(\eta(n))]-k\right)\psi(k)\cos
kt=
$$
$$
=\frac{1}{[\eta(n)]-n}\sum\limits_{k=n+1}^{[\eta(n)]-1}(k-n)\psi(k)\cos
kt+ \psi([\eta(n)])\cos([\eta(n)]t)+
$$
\begin{equation}\label{ww}
+\frac{1}{[\eta(\eta(n))]-[\eta(n)]}
\sum\limits_{k=[\eta(n)]+1}^{[\eta(\eta(n))]-1}\left([\eta(\eta(n))]-k\right)\psi(k)\cos
kt.
\end{equation}

Оскільки $W_{N,M}(\lambda;\gamma;t)$ є тригонометричним поліномом
порядку $M$, то згідно з означенням $(\psi,\beta)$--похідної, серед
функцій  $\left(W_{N,M}(\lambda;\gamma;t)\right)_{\beta}^{\psi}$
знайдеться така, що також є тригонометричним поліномом порядку $M$.
Надалі саме таку функцію будемо розуміти під записом
$\left(W_{N,M}(\lambda;\gamma;t)\right)_{\beta}^{\psi}$.

Враховуючи співвідношення (\ref{ww}) та означення
$(\psi,\beta)$--похідної, для довільних $\psi\in\mathfrak{M}$,
$\beta\in\mathbb{R}$  отримуємо рівність
$$
\left(W_{[\eta(n)],[\eta(\eta(n))]}(\psi;0;t)-
W_{n,[\eta(n)]}(\psi;0;t)\right)^{\psi}_{\beta}=
$$
$$
=\frac{1}{[\eta(n)]-n}\sum\limits_{k=n+1}^{[\eta(n)]-1}(k-n)\cos\!\left(kt+\frac{\beta\pi}{2}\right)
+\cos\left([\eta(n)]t+\frac{\beta\pi}{2}\right)+
$$
\begin{equation}\label{1t5}
+\frac{1}{[\eta(\eta(n))]-[\eta(n)]}
\sum\limits_{k=[\eta(n)]+1}^{[\eta(\eta(n))]-1}\left([\eta(\eta(n))]-k\right)\cos
\left(kt+\frac{\beta\pi}{2}\right).
\end{equation}

Для доведення нерівності (\ref{gg11}) буде корисним наступне
твердження.

\textbf{Лема 2.} \emph{Нехай $\psi\in\mathfrak{M}^{+}_{\infty}$,
$\eta(n)-n\geq a$, $\mu(n)\geq b$, $\beta$ --- довільне дійсне
число. Тоді }

\noindent\emph{ 1) якщо $a>0$, $b>0$, то}
$$
\left|\left(W_{[\eta(n)],[\eta(\eta(n))]}(\psi;0;t)-
W_{n,[\eta(n)]}(\psi;0;t)\right)^{\psi}_{\beta}\right|\leq
$$
\begin{equation}\label{incl1}
\leq \left(1+\frac{1}{2a}+\frac{1}{2b}\right)(\eta(n)-n), \
t\in\mathbb{R};
\end{equation}
\noindent\emph{ 2) якщо $a>2$, $b>0$, то}
$$
\left|\left(W_{[\eta(n)],[\eta(\eta(n))]}(\psi;0;t)-
W_{n,[\eta(n)]}(\psi;0;t)\right)^{\psi}_{\beta}\right|\leq
$$
\begin{equation}\label{incl2}
\leq
\left(\frac{a}{a-1}+\frac{2a}{a-2}\right)\frac{\pi^{2}}{t^{2}}\frac{1}{\eta(n)-n},
 \ \  0<|t|\leq\pi.
\end{equation}

 \emph{Доведення леми 2.} Спочатку встановимо істинність оцінки
 (\ref{incl1}). З рівності (\ref{1t5}), випливає, що
$$
\left|\left(W_{[\eta(n)],[\eta(\eta(n))]}(\psi;0;t)\!-\!
W_{n,[\eta(n)]}(\psi;0;t)\right)^{\psi}_{\beta}\right|\!\leq\!
\frac{1}{[\eta(n)]\!-\!n}\!\!\sum\limits_{k=n+1}^{[\eta(n)]-1}\!\!\!(k-n)+
$$
$$
+1+\frac{1}{[\eta(\eta(n))]-[\eta(n)]}
\sum\limits_{k=[\eta(n)]+1}^{[\eta(\eta(n))]-1}\left([\eta(\eta(n))]-k\right)
=\frac{[\eta(n)]-n-1}{2}+1+
$$
\begin{equation}\label{gg0}
+\frac{[\eta(\eta(n))]-[\eta(n)]-1}{2}
=\frac{1}{2}\left(([\eta(\eta(n))]-[\eta(n)])+([\eta(n)]-n)\right).
\end{equation}

Щоб оцінити величини $[\eta(\eta(n))]-[\eta(n)]$ i $[\eta(n)]-n$
використаємо наступне твердження.

\textbf{Лема 3.} \emph{Нехай $\psi\in\mathfrak{M}^{+}_{\infty}$,
$\eta(n)-n\geq a>0$, $\mu(n)\geq b>0$. Тоді}

\noindent\emph{ 1) якщо $a>1$, $b>0$, то}
\begin{equation}\label{i2}
\left(1-\frac{1}{a}\right)\left(\eta(n)-n\right)<[\eta(n)]-n\leq
\eta(n)-n;
\end{equation}
\emph{ 2) якщо $a>2$, $b>0$, то}
\begin{equation}\label{i3}
\left(\frac{1}{2}-\frac{1}{a}\right)\left(\eta(n)-n\right)<[\eta(\eta(n))]-[\eta(n)]<
\left(1+\frac{1}{a}+\frac{1}{b}\right)\left(\eta(n)-n\right).
\end{equation}

\emph{Доведення леми 3.} Друга нерівність в (\ref{i2}) є очевидною.
Позначивши  $\{\alpha\}=\alpha-[\alpha]$, при ${\eta(n)-n\geq a>1}$
одержуємо
$$
[\eta(n)]-n=\eta(n)-n-\left\{\eta(n)\right\}=\left(\eta(n)-n\right)\left(1-\frac{\left\{\eta(n)\right\}}{\eta(n)-n}\right)>
$$
$$
>\left(1-\frac{1}{a}\right)\left(\eta(n)-n\right).
$$

Тим самим (\ref{i2}) доведено.

 Щоб переконатись в справедливості нерівностей (\ref{i3}), спочатку покажемо, що
при $ \eta(n)-n\geq a>0$ і $\mu(n)\geq b>0$ має місце співвідношення
\begin{equation}\label{i1}
\frac{1}{2}\left(\eta(n)-n\right)\leq
\eta(\eta(n))-\eta(n)<\left(1+\frac{1}{b}\right)\left(\eta(n)-n\right).
\end{equation}

Дійсно, беручи до уваги означення функції $\mu(t)$, для довільної
$\psi\in\mathfrak{M}$ справедлива рівність
\begin{equation}\label{gg1}
\eta(t)=t\left(1+\frac{\eta(t)-t}{t}\right)=t\left(1+\frac{1}{\mu(t)}\right).
\end{equation}

Оскільки $\psi\in\mathfrak{M}^{+}_{\infty}$, то функція
$\frac{1}{\mu(t)}$ монотонно прямує до нуля при
$t\rightarrow\infty$. Нехай $\mu(t)\geq b>0$. Тоді в  силу
(\ref{gg1}), маємо
\begin{equation}\label{in1}
\eta'(t)=1+\frac{1}{\mu(t)}+t\left(\frac{1}{\mu(t)}\right)'\leq
1+\frac{1}{\mu(t)}\leq 1+\frac{1}{b}.
\end{equation}

Зазначимо також, що для довільної функції $\psi\in\mathfrak{M}$
(див., наприклад, \cite[c. 162--163]{Stepanets1})
\begin{equation}\label{in2}
\eta'(t)=\frac{\psi'(t)}{2\psi'(\eta(t))}\geq \frac{1}{2}, \ t\geq
1, \ \eta'(t)=\eta'(t+0).
\end{equation}
З  (\ref{in1}) i (\ref{in2}), а також з рівності
$$
\eta(\eta(t))-\eta(t)=\int\limits_{t}^{\eta(t)}\eta'(u)du,
$$
випливає (\ref{i1}).

Застосування  нерівностей (\ref{i1}) при ${a>0, \ b>0}$ дозволяє
записати співвідношення
$$
[\eta(\eta(n))]-[\eta(n)]\leq
\eta(\eta(n))-\eta(n)+\left\{\eta(n)\right\}<
\left(1+\frac{1}{b}\right)\left(\eta(n)-n\right)+1=
$$
\begin{equation}\label{gg2}
=\left(\eta(n)-n\right)\left(1+\frac{1}{b}+\frac{1}{\eta(n)-n}\right)
\leq\left(1+\frac{1}{a}+\frac{1}{b}\right)\left(\eta(n)-n\right),
\end{equation}
а при $a>2, \ b>0$ --- співвідношення
$$
[\eta(\eta(n))]-[\eta(n)]\geq
\eta(\eta(n))-\eta(n)-\left\{\eta(\eta(n))\right\}>
\frac{1}{2}\left(\eta(n)-n\right)-1=
$$
\begin{equation}\label{gg3}
= \left(\eta(n)-n\right)\left(\frac{1}{2}-\frac{1}{\eta(n)-n}\right)
\geq\left(\frac{1}{2}-\frac{1}{a}\right)\left(\eta(n)-n\right).
\end{equation}

Нерівності (\ref{gg2}) i (\ref{gg3}) доводять (\ref{i3}). Лему 3
доведено.

Щоб переконатись у справедливості (\ref{incl1}), достатньо
скористатись формулою (\ref{gg0}) і застосувати нерівності
(\ref{i2}) і (\ref{i3}) леми 3.

Перейдемо до доведення нерівності (\ref{incl2}). В силу означення
$(\psi,\beta)$--похідної та на підставі (\ref{www}) для довільної
$\psi\in\mathfrak{M}$  одержуємо
$$
\left(W_{[\eta(n)],[\eta(\eta(n))]}(\psi;0;t)-
W_{n,[\eta(n)]}(\psi;0;t)\right)^{\psi}_{\beta}=
$$
$$
=\left(\frac{1}{[\eta(\eta(n))]-[\eta(n)]}
\sum\limits_{k=[\eta(n)]}^{[\eta(\eta(n))]-1}
\sum\limits_{j=1}^{k}\psi(j)\cos jt\right)^{\psi}_{\beta}-
$$
$$
-\left(\frac{1}{[\eta(n)]-n}\sum\limits_{k=n}^{[\eta(n)]-1}\sum\limits_{j=1}^{k}\psi(j)\cos
jt\right)^{\psi}_{\beta}=
$$
$$
=\frac{1}{[\eta(\eta(n))]-[\eta(n)]}
\sum\limits_{k=[\eta(n)]}^{[\eta(\eta(n))]-1}\sum\limits_{j=1}^{k}\cos\left(jt+\frac{\beta\pi}{2}\right)-
$$
\begin{equation}\label{t2}
-\frac{1}{[\eta(n)]-n}\sum\limits_{k=n}^{[\eta(n)]-1}
\sum\limits_{j=1}^{k}\cos\left(jt+\frac{\beta\pi}{2}\right).
\end{equation}

Застосовуючи до правої  частини рівності (\ref{t2}) рівність
(\ref{gg8}), та формулу
$$
\sum\limits_{k=0}^{N-1}\sin\left(x+ky\right)=\sin\left(x+\frac{N-1}{2}y\right)\sin
Ny\cosec \frac{y}{2},
$$
(див., наприклад, \cite[с. 43]{Gradshteyn}) із (\ref{t2}) отримуємо
$$
\left(W_{[\eta(n)],[\eta(\eta(n))]}(\psi;0;t)-
W_{n,[\eta(n)]}(\psi;0;t)\right)^{\psi}_{\beta}=
$$
$$
=\frac{1}{[\eta(\eta(n))]-[\eta(n)]}
\sum\limits_{k=[\eta(n)]}^{[\eta(\eta(n))]-1}D_{k,-\beta}(t)-\frac{1}{[\eta(n)]-n}\sum\limits_{k=n}^{[\eta(n)]-1}D_{k,-\beta}(t)=
$$
$$
=\frac{1}{[\eta(\eta(n))]-[\eta(n)]}
\sum\limits_{k=[\eta(n)]}^{[\eta(\eta(n))]-1}\frac{\sin\left(\left(k+\frac{1}{2}\right)t+
\frac{\beta\pi}{2}\right)-\sin\frac{\beta\pi}{2}\cos\frac{t}{2}}{2\sin\frac{t}{2}}-
$$
$$
-\frac{1}{[\eta(n)]-n}\sum\limits_{k=n}^{[\eta(n)]-1}\frac{\sin\left(\left(k+\frac{1}{2}\right)t+
\frac{\beta\pi}{2}\right)-\sin\frac{\beta\pi}{2}\cos\frac{t}{2}}{2\sin\frac{t}{2}}=
$$
$$ =\frac{1}{2\sin\frac{t}{2}}\left(\frac{1}{[\eta(\eta(n))]-[\eta(n)]}
\sum\limits_{k=[\eta(n)]}^{[\eta(\eta(n))]-1}\sin\left(\left(k+\frac{1}{2}\right)t+
\frac{\beta\pi}{2}\right)-\right.
$$
$$
\left.-\frac{1}{[\eta(n)]-n}\sum\limits_{k=n}^{[\eta(n)]-1}\sin\left(\left(k+\frac{1}{2}\right)t+
\frac{\beta\pi}{2}\right)\right)=
$$
$$=\!\frac{1}{2\left(\sin\frac{t}{2}\right)^{2}}\!\left(\!\frac{1}{[\eta(\eta(n))]\!-\![\eta(n)]}
\!\left(\!\sin\!\left(\frac{[\eta(\eta(n))]}{2}t\!+\!
\frac{\beta\pi}{2}\right)\!\sin\frac{[\eta(\eta(n))]}{2}t-
\right.\right.
$$
$$
\left.-\sin\left(\frac{[\eta(n)]}{2}t+
\frac{\beta\pi}{2}\right)\sin\frac{[\eta(n)]}{2}t \right)-
$$
$$ -\frac{1}{[\eta(n)]-n}
\left(\sin\left(\frac{[\eta(n)]}{2}t+\frac{\beta\pi}{2}\right)\sin\frac{[\eta(n)]}{2}t-
\right.
$$
\begin{equation}\label{t8}
\left.\left.-\sin\left(\frac{n}{2}t+
\frac{\beta\pi}{2}\right)\sin\frac{n}{2}t \right)\right).
\end{equation}

Беручи до уваги рівність (\ref{t8}) i нерівність (\ref{ineq}), маємо
$$
\left|\left(W_{[\eta(n)],[\eta(\eta(n))]}(\psi;0;t)-
W_{n,[\eta(n)]}(\psi;0;t)\right)^{\psi}_{\beta}\right|\leq
$$
\begin{equation}\label{t3}
\leq\frac{\pi^{2}}{t^{2}}\left(\frac{1}{[\eta(\eta(n))]-[\eta(n)]}+\frac{1}{[\eta(n)]-n}\right),
\ \  0<|t|\leq \pi.
\end{equation}
Враховуючи (\ref{t3}), та використовуючи нерівності (\ref{i2}) i
(\ref{i3}), отримуємо (\ref{incl2}). Лему 2 доведено.

Використовуючи лему 2, покажемо виконання нерівності (\ref{gg11}). З
(\ref{incl1}) випливає оцінка
$$
\int\limits_{|t|\leq
\frac{1}{\eta(n)-n}}\left|\left(W_{[\eta(n)],[\eta(\eta(n))]}(\psi;0;t)-
W_{n,[\eta(n)]}(\psi;0;t)\right)^{\psi}_{\beta}\right|^{p}dt\leq
$$
$$
\leq \int\limits_{|t|\leq
\frac{1}{\eta(n)-n}}\left(1+\frac{1}{2a}+\frac{1}{2b}\right)^{p}(\eta(n)-n)^{p}dt=
$$
\begin{equation}\label{t6}
= 2\left(1+\frac{1}{2a}+\frac{1}{2b}\right)^{p}(\eta(n)-n)^{p-1},
\end{equation}
а з нерівності (\ref{incl2}) --- оцінка
$$
\int\limits_{\frac{1}{\eta(n)-n}\leq |t|\leq\pi
}\left|\left(W_{[\eta(n)],[\eta(\eta(n))]}(\psi;0;t)-
W_{n,[\eta(n)]}(\psi;0;t)\right)^{\psi}_{\beta}\right|^{p}dt\leq
$$
$$
\leq \left(\frac{a}{a-1}+\frac{2a}{a-2}\right)^{p}\pi^{2p}
\frac{1}{\left(\eta(n)-n\right)^{p}}\int\limits_{
\frac{1}{\eta(n)-n}\leq |t|\leq \pi}\frac{1}{t^{2p}}dt=
$$
$$
=\left(\frac{a}{a\!-\!1}\!+\!\frac{2a}{a\!-\!2}\right)^{p}\pi^{2p}
(\eta(n)-n)^{p-1}\frac{2}{2p-1}\left(1-
\frac{\pi^{1-2p}}{(\eta(n)-n)^{2p-1}}\right)<
$$
\begin{equation}\label{t7}
<
2\pi^{2p}\left(\frac{a}{a-1}+\frac{2a}{a-2}\right)^{p}(\eta(n)-n)^{p-1}.
\end{equation}
Об'єднуючи (\ref{t6})--(\ref{t7}), та враховуючи очевидну нерівність
\begin{equation}\label{t11}
\frac{a}{a-1}+\frac{2a}{a-2}>1+\frac{1}{2a}+\frac{1}{2b}, \ a>2, \
b>2,
\end{equation}
маємо
$$
\left\|\left(W_{[\eta(n)],[\eta(\eta(n))]}(\psi;0;t)-
W_{n,[\eta(n)]}(\psi;0;t)\right)^{\psi}_{\beta}\right\|_{p}\leq
$$
$$
\leq\left(\frac{a}{a-1}+\frac{2a}{a-2}\right)2^{\frac{1}{p}}
\left(1+\pi^{2p}\right)^{\frac{1}{p}}(\eta(n)-n)^{1-\frac{1}{p}}\leq
$$
$$
\leq2\left(\frac{a}{a-1}+\frac{2a}{a-2}\right)
\left(1+\pi^{2}\right)(\eta(n)-n)^{1-\frac{1}{p}}=
$$
\begin{equation}\label{1t6}
=\frac{2\left(1+\pi^{2}\right)a(3a-4)}{(a-1)(a-2)}(\eta(n)-n)^{1-\frac{1}{p}},
\ 1\leq p<\infty.
\end{equation}

Крім того, з (\ref{incl1}) та (\ref{t11}) випливає, що при $a>2$ i
$b>2$
$$
\left\|\left(W_{[\eta(n)],[\eta(\eta(n))]}(\psi;0;t)-
W_{n,[\eta(n)]}(\psi;0;t)\right)^{\psi}_{\beta}\right\|_{\infty}\leq
$$
\begin{equation}\label{t9}
\leq\left(1+\frac{1}{2a}+\frac{1}{2b}\right)(\eta(n)-n)<
\frac{a(3a-4)}{(a-1)(a-2)}(\eta(n)-n).
\end{equation}

 З
(\ref{function}),  (\ref{1t6}) i (\ref{t9}) одержуємо оцінку
$$
\Big\|\big(f_{p}(t)\big)^{\psi}_{\beta}\Big\|_{p}=
\frac{(a-1)(a-2)}{2\left(1+\pi^{2}\right)a(3a-4)}\frac{1}{(\eta(n)-n)^{1-\frac{1}{p}}}\times
$$
\begin{equation}\label{norm}
\times\Big\|\left(W_{[\eta(n)],[\eta(\eta(n))]}(\psi;0;t)-
W_{n,[\eta(n)]}(\psi;0;t)\right)^{\psi}_{\beta}\!\Big\|_{p} \leq 1,
\ 1\leq p\leq \infty,
\end{equation}
яка доводить  включення $f_{p}\in C^{\psi}_{\beta,p},$ при всіх $
1\leq p\leq \infty$ і $\beta\in\mathbb{R}$.

Далі покажемо, що шуканою функцією $f^{*}$ є функція $f_{p}$, тобто
для $f^{*}(\cdot)=f_{p}(\cdot)$ виконується оцінка (\ref{t10}).
Дійсно, в силу (\ref{ww}), бачимо, що для будь--якого
тригонометричного полінома
 $t_{n-1}\in~ \mathcal{T}_{2n-1}$
\begin{equation}\label{t4}
\int\limits_{-\pi}^{\pi}\left(W_{[\eta(n)],[\eta(\eta(n))]}(1;0;t)-
W_{n,[\eta(n)]}(1;0;t)\right) t_{n-1}(t)dt=0.
\end{equation}

Із (\ref{function}) та (\ref{t4}) маємо
$$
\int\limits_{-\pi}^{\pi}\big(f_{p}(t)-t_{n-1}(t)\big)\left(W_{[\eta(n)],[\eta(\eta(n))]}(1;0;t)-
W_{n,[\eta(n)]}(1;0;t)\right)dt =
$$
$$
=\int\limits_{-\pi}^{\pi}f_{p}(t)\left(W_{[\eta(n)],[\eta(\eta(n))]}(1;0;t)-
W_{n,[\eta(n)]}(1;0;t)\right)dt=
$$
$$
=\frac{(a-1)(a-2)}{2\left(1+\pi^{2}\right)a(3a-4)}\frac{1}{(\eta(n)-n)^{1-\frac{1}{p}}}\times
$$
$$
\times\int\limits_{-\pi}^{\pi}\left(W_{[\eta(n)],[\eta(\eta(n))]}(\psi;0;t)-
W_{n,[\eta(n)]}(\psi;0;t)\right)\times
$$
\begin{equation}\label{yyy}
\times\left(W_{[\eta(n)],[\eta(\eta(n))]}(1;0;t)-
W_{n,[\eta(n)]}(1;0;t)\right)dt.
\end{equation}

Застосовуючи рівність (\ref{lemma1}) до функцій
$W_{N,M}(\lambda;\gamma;t)$  при ${\lambda(k)=1}, \ \gamma=0$,
$N=n$, $M=[\eta(n)]$, а також при $\lambda(k)=1, \ \gamma=0$,
$N=[\eta(n)]$, $M=[\eta(\eta(n))]$, і діючи так само, як і при
доведенні співвідношення (\ref{ww}), легко показати, що
$$
W_{[\eta(n)],[\eta(\eta(n))]}(1;0;t)- W_{n,[\eta(n)]}(1;0;t)=
$$
$$
=\frac{1}{[\eta(n)]-n}\sum\limits_{k=n+1}^{[\eta(n)]-1}(k-n)\cos kt+
\cos([\eta(n)]t)+
$$
\begin{equation}\label{wwww}
+\frac{1}{[\eta(\eta(n))]-[\eta(n)]}
\sum\limits_{k=[\eta(n)]+1}^{[\eta(\eta(n))]-1}\left([\eta(\eta(n))]-k\right)\cos
kt.
\end{equation}
На основі рівностей (\ref{ww}) і (\ref{wwww}), а також формул
$$
\int\limits_{-\pi}^{\pi}\cos kt\cos mtdt={\left\{\begin{array}{cc}
0, \ & k\neq m, \\
\pi, & k=m, \
  \end{array} \right.} \ \ k,m\in\mathbb{N},
$$
отримуємо
$$
\int\limits_{-\pi}^{\pi}\left(W_{[\eta(n)],[\eta(\eta(n))]}(\psi;0;t)-
W_{n,[\eta(n)]}(\psi;0;t)\right)\times
$$
$$
\times\left(W_{[\eta(n)],[\eta(\eta(n))]}(1;0;t)-
W_{n,[\eta(n)]}(1;0;t)\right)dt=
$$
$$
=\frac{\pi}{([\eta(n)]-n)^{2}}\!\!\!\sum\limits_{k=n+1}^{[\eta(n)]-1}\psi(k)(k-n)^{2}
\!+
$$
\begin{equation}\label{gg5}
+\pi\psi([\eta(n)])+\frac{\pi}{([\eta(\eta(n))]-[\eta(n)])^{2}}
\sum\limits_{k=[\eta(n)]+1}^{[\eta(\eta(n))]-1}\psi(k)\left([\eta(\eta(n))]-k\right)^{2}.
\end{equation}

Оскільки $\psi(t)$ монотонно спадає, то, згідно з означенням
характеристики $\eta(t)$, маємо
$$
\frac{\pi}{([\eta(n)]-n)^{2}}\sum\limits_{k=n+1}^{[\eta(n)]-1}\psi(k)(k-n)^{2}
+\pi\psi([\eta(n)])+
$$
$$+\frac{\pi}{([\eta(\eta(n))]-[\eta(n)])^{2}}
\sum\limits_{k=[\eta(n)]+1}^{[\eta(\eta(n))]-1}\psi(k)\left([\eta(\eta(n))]-k\right)^{2}\geq
$$
$$
>
\pi\psi(\eta(\eta(n))])\left(\frac{1}{([\eta(n)]-n)^{2}}\sum\limits_{k=n+1}^{[\eta(n)]-1}(k-n)^{2}
+1+\right.
$$
$$\left.+\frac{1}{([\eta(\eta(n))]-[\eta(n)])^{2}}
\sum\limits_{k=[\eta(n)]+1}^{[\eta(\eta(n))]-1}\left([\eta(\eta(n))]-k\right)^{2}\right)=
$$
$$
=
\frac{\pi}{4}\psi(n)\left(\frac{1}{([\eta(n)]-n)^{2}}\sum\limits_{k=1}^{[\eta(n)]-n-1}k^{2}+1+\right.
$$
\begin{equation}\label{t5}
\left.+\frac{1}{([\eta(\eta(n))]\!-\![\eta(n)])^{2}}
\sum\limits_{k=1}^{[\eta(\eta(n))]-[\eta(n)]-1}k^{2}\right).
\end{equation}

Використовуючи формулу (див., наприклад, \cite[c. 15]{Gradshteyn})
$$
\sum\limits_{k=1}^{M}k^{2}=\frac{M(M+1)(2M+1)}{6}, \ \
M\in\mathbb{N},
$$
при $M=[\eta(n)]-n-1$ та  $M=[\eta(\eta(n))]-[\eta(n)]-1$,
одержуємо
$$
\frac{1}{([\eta(n)]-n)^{2}}\sum\limits_{k=1}^{[\eta(n)]-n-1}k^{2}+1+
\frac{1}{([\eta(\eta(n))]\!-\![\eta(n)])^{2}}
\sum\limits_{k=1}^{[\eta(\eta(n))]-[\eta(n)]-1}k^{2}=
$$
$$
=\frac{([\eta(n)]-n-1)([\eta(n)]-n)(2[\eta(n)]-2n-1)}{6([\eta(n)]-n)^{2}}+1+
$$
$$
+\frac{([\eta(\eta(n))]\!-\![\eta(n)]-1)([\eta(\eta(n))]\!-\![\eta(n)])(2[\eta(\eta(n))]\!-\!2[\eta(n)]\!-\!1)}
{6([\eta(\eta(n))]-[\eta(n)])^{2}}=
$$
$$
=\frac{([\eta(n)]-n-1)(2[\eta(n)]-2n-1)}{6([\eta(n)]-n)}+1 +
$$
$$
+\frac{([\eta(\eta(n))]-[\eta(n)]-1)(2[\eta(\eta(n))]-2[\eta(n)]-1)}{6([\eta(\eta(n))]-[\eta(n)])}=
$$
$$
=\frac{1}{6}\left(2([\eta(n)]-n)+2([\eta(\eta(n))]-[\eta(n)])+\frac{1}{[\eta(n)]-n}+\right.
$$
\begin{equation}\label{t12}
\left.+\frac{1}{[\eta(\eta(n))]-[\eta(n)]}\right).
\end{equation}

В силу нерівностей (\ref{i2}) i (\ref{i3}), при $a>2$, маємо
$$
\frac{1}{6}\left(2([\eta(n)]-n)+2([\eta(\eta(n))]-[\eta(n)])+\frac{1}{[\eta(n)]-n}+\right.
$$
\begin{equation}\label{gg6}
\left.+\frac{1}{[\eta(\eta(n))]-[\eta(n)]}\right)>\frac{3a-4}{6a}(\eta(n)-n).
\end{equation}

Об'єднання формул (\ref{yyy})--(\ref{gg6}) дозволяє записати
наступну оцінку, справедливу для будь-якого
$t_{n-1}\in\mathcal{T}_{2n-1}$:
$$
\int\limits_{-\pi}^{\pi}\big(f_{p}(t)-t_{n-1}(t)\big)\left(W_{[\eta(n)],[\eta(\eta(n))]}(1;0;t)-
W_{n,[\eta(n)]}(1;0;t)\right)dt\geq
$$
\begin{equation}\label{t14}
\geq\frac{\pi(a-1)(a-2)}{48\left(1+\pi^{2}\right)a^{2}}
\psi(n)(\eta(n)-n)^{\frac{1}{p}}.
\end{equation}

 З іншого боку, в силу (\ref{www}) і означення $(\psi,\beta)$--похідної при $\beta=0$
  $$W_{[\eta(n)],[\eta(\eta(n))]}(1;0;t)-
W_{n,[\eta(n)]}(1;0;t)=$$
 \begin{equation}\label{riv}
=\left(W_{[\eta(n)],[\eta(\eta(n))]}(\psi;0;t)-
W_{n,[\eta(n)]}(\psi;0;t)\right)^{\psi}_{0},
\end{equation}
тому використовуючи нерівності (\ref{statement1})  i (\ref{1t6}),
отримуємо
$$
\int\limits_{-\pi}^{\pi}\big(f_{p}(t)-t_{n-1}(t)\big)\left(W_{[\eta(n)],[\eta(\eta(n))]}(1;0;t)-
W_{n,[\eta(n)]}(1;0;t)\right)dt \leq
$$
$$
\leq\|f_{p}(t)-t_{n-1}(t)\|_{\infty}\left\|\left(W_{[\eta(n)],[\eta(\eta(n))]}(1;0;t)-
W_{n,[\eta(n)]}(1;0;t)\right)\right\|_{1}\leq
$$
\begin{equation}\label{gg7}
\leq
\frac{2\left(1+\pi^{2}\right)a(3a-4)}{(a-1)(a-2)}\|f_{p}(t)-t_{n-1}(t)\|_{\infty}.
\end{equation}

З (\ref{t14}) i (\ref{gg7}) випливає, що для довільного
$t_{n-1}\in\mathcal{T}_{2n-1}$
$$
\|f_{p}(t)-t_{n-1}(t)\|_{\infty}
\geq\frac{\pi(a-1)^{2}(a-2)^{2}}{96\left(1+\pi^{2}\right)^{2}a^{3}(3a-4)}
\psi(n)(\eta(n)-n)^{\frac{1}{p}}=
$$
\begin{equation}\label{zz3}
=C_{a}\psi(n)(\eta(n)-n)^{\frac{1}{p}}.
\end{equation}

Із (\ref{zz3}) для класів $C^{\psi}_{\beta,p}$ при
$\psi\in\mathfrak{M}^{+}_{\infty}, \
\lim\limits_{t\rightarrow\infty}(\eta(\psi,t)-t)=\infty$,
${\eta(n)-n\geq a>2}, \ {\mu(n)\geq b>2}$, $\beta\in \mathbb{R}$
одержуємо оцінку (\ref{t10}). Теорему 1 доведено.

Важливим прикладом функцій $\psi(t)$ з множини
$\mathfrak{M}^{+}_{\infty}$, які задовольняють умову
$\lim\limits_{t\rightarrow\infty}(\eta(\psi,t)-t)=\infty$, є функції
\begin{equation}\label{psi_r}
{\psi_{r}(t)=\exp(-\alpha t^{r})}, \ {\alpha>0}, \ {r\in(0,1)}.
\end{equation}

Для них $\eta(\psi_{r};n)=\left(\alpha^{-1}\ln
2+n^{r}\right)^{\frac{1}{r}}$. Тоді, використавши узагальнену
нерівність Бернуллі
$$
(1+x)^{\rho}\geq1+\rho x, \ x>-1, \
\rho\in(-\infty,0]\cup[1,\infty),
$$
отримуємо
\begin{equation}\label{eta-n}
\eta(\psi_{r};n)-n=
 n\left(\left(1+\frac{\ln 2}{\alpha n^{r}}\right)^{\frac{1}{r}}-1\right)\geq \frac{\ln 2}{\alpha r}n^{1-r}, \ n\in\mathbb{N}.
\end{equation}

 З формули (\ref{eta-n}) випливає, що  для всіх номерів $n\geq1+\left(\frac{2r\alpha}{\ln2}\right)^{\frac{1}{1-r}}$ виконується нерівність
$$
\eta(\psi_{r};n)-n\geq a>2,
$$
при
\begin{equation}\label{aa}
a=a(\alpha,r)=\frac{\ln 2}{\alpha
r}\left(1+\left(\frac{2r\alpha}{\ln2}\right)^{\frac{1}{1-r}}\right)^{1-r}.
\end{equation}

В силу (\ref{aa})
 $$\mu(\psi_{r};n)=\frac{n}{\eta(\psi_{r};n)-n}=\frac{1}{\left(\frac{\ln 2}{\alpha n^{r}}+1\right)^{\frac{1}{r}}-1}$$
і, як неважко переконатись, для всіх $n\geq
1+2\left(\frac{\ln2}{\alpha\left(3^{r}-2^{r}\right)}\right)^{\frac{1}{r}}$
 виконується нерівність
$$
\mu(\psi_{r};n)\geq b>2,
$$
де
\begin{equation}\label{bb}
b=b(r,\alpha)=\left(\left(\frac{\ln 2}{\alpha
}\left(1+2\left(\frac{\ln2}{\alpha\left(3^{r}-2^{r}\right)}\right)
^{\frac{1}{r}}\right)^{-r}+1\right)^{\frac{1}{r}}-1\right)^{-1}.
\end{equation}

З наведених вище міркувань випливає, що до класів
$C^{\psi}_{\beta,p}$, породжених послідовностями $\psi_{r}(t)$
вигляду (\ref{psi_r}) можна застосувати теорему 1, в умові якої
параметри $a$ i $b$ визначаються формулами (\ref{aa}) i (\ref{bb})
відповідно. В результаті одержимо наступне твердження.

\textbf{Наслідок 1.} \emph{ Нехай $\psi_{r}(t)=\exp\left(-\alpha
t^{r}\right), \  {r\in(0,1)}, \ {\alpha>0}$, ${1\leq  p< \infty}$,
$\beta\in \mathbb{R}$. Тоді  для усіх номерів $n$, таких, що
$$n\geq1+\max\left\{\left(\frac{2r\alpha}{\ln2}\right)^{\frac{1}{1-r}},
2\left(\frac{\ln2}{\alpha\left(3^{r}-2^{r}\right)}\right)^{\frac{1}{r}}\right\},$$
 справедливі оцінки}
$$
C_{a} \exp\left(-\alpha
n^{r}\right)n^{\frac{1}{p}}\left(\left(\frac{\ln2}{\alpha
n^{r}}+1\right)^{\frac{1}{r}}-1\right)^{\frac{1}{p}} \leq \ \ \ \ \
\ \ \ \ \ \ \ \ \ \ \ \ \ \ \ \ \ \ \ \ \ \ \ \ \ \ \
$$
$$
\leq{ E}_{n}\Big(C^{\psi_{r}}_{\beta,p}\Big)_{C}\leq{\cal
E}_{n}\big(C^{\psi_{r}}_{\beta,p})_{C}\leq
$$
\begin{equation}\label{cons}
\ \ \ \ \ \ \ \ \ \ \ \ \ \ \ \ \ \ \ \ \ \ \ \ \ \ \ \ \ \  \leq
C_{a,b} \ (2p)^{1-\frac{1}{p}} \exp\left(-\alpha
n^{r}\right)n^{\frac{1}{p}}\left(\left(\frac{\ln2}{\alpha
n^{r}}+1\right)^{\frac{1}{r}}-1\right)^{\frac{1}{p}},
 \end{equation}
\emph {де величини $C_{a}$ i $C_{a,b}$ означаються формулами
(\ref{Ca}) i (\ref{Cab}) при ${a=a(\alpha,r)}$, $b=b(\alpha,r)$, що
задані за допомогою рівностей (\ref{aa}) i (\ref{bb}) відповідно.}

Зазначимо також, що оскільки
$$\eta(\psi_{r};n)-n=n\left(\left(\frac{\ln2}{\alpha
n^{r}}+1\right)^{\frac{1}{r}}-1\right)\asymp n^{1-r}, \ {r\in(0,1]},
\ {\alpha>0},$$
 то з (\ref{cons}) випливають порядкові рівності
\begin{equation}\label{consequence}
\mathcal{E}_{n}\Big(C^{\psi_{r}}_{\beta,p}\Big)_{C}\asymp{
E}_{n}\Big(C^{\psi_{r}}_{\beta,p}\Big)_{C}\asymp\exp\left(-\alpha
n^{r}\right)n^{\frac{1-r}{p}},\ 1\leq p< \infty,
 \end{equation}
де запис $A(n)\asymp B(n)$ $(A(n)>0, \ B(n)>0)$ означає існування
додатних сталих $K_{1}$ і $K_{2}$ таких, що ${K_{1}B(n)\leq A(n)\leq
K_{2}B(n)}$, $n\in\mathbb{N}$.

Для величини $\mathcal{E}_{n}\Big(C^{\psi_{r}}_{\beta,p}\Big)_{C}$
оцінка (\ref{consequence}) встановлена в роботі \cite{Rom}.

\textbf{Теорема 2.}\emph{ Нехай $\psi\in\mathfrak{M}^{+}_{\infty}, \
\lim\limits_{t\rightarrow\infty}(\eta(\psi,t)-t)=\infty$, ${1<
s\leq\infty}$, ${\frac{1}{s}+\frac{1}{s'}=1}$, $\beta\in
\mathbb{R}$. Тоді  для довільних $n\in
 \mathbb{N}$, таких, що $ {\eta(n)-n\geq a>2}, \ {\mu(n)\geq b>2}$ справедливі оцінки}
$$
C_{a}\psi(n)(\eta(n)-n)^{\frac{1}{s'}} \leq
E_{n}\big(L^{\psi}_{\beta,1}\big)_{s}\leq{\cal
E}_{n}\big(L^{\psi}_{\beta,1})_{s} \leq \ \ \ \ \ \ \ \ \ \ \ \ \ \
\ \ \ \ \ \ \ \ \ \ \ \ \ \ \ \ \ \ \ \
$$
\begin{equation}\label{theorem2}
\ \ \ \ \ \ \ \ \ \ \ \ \ \ \ \ \ \ \ \ \ \ \ \ \ \ \ \ \ \ \ \ \ \
\ \ \ \  \ \ \ \ \ \ \ \ \ \ \ \ \ \ \ \ \ \ \ \ \leq C_{a,b} \
\left(2s'\right)^{\frac{1}{s}}\psi(n)\left(\eta(n)-n\right)^{\frac{1}{s'}},
\end{equation}
\emph{ де величини $C_{a}$ i $C_{a,b}$ означаються формулами
(\ref{Ca}) i (\ref{Cab}) відповідно.}

\textbf{\emph{Доведення.}} Для отримання оцінки зверху величини
${\cal E}_{n}\big(L^{\psi}_{\beta,1})_{s}$ використаємо інтегральне
зображення (\ref{2t1}) та нерівність Юнга (див., наприклад, \cite[c.
293]{Stepanets1}). Тоді для довільних $1\leq s\leq\infty$, одержимо
\begin{equation}\label{3t1}
{\cal E}_{n}\big(L^{\psi}_{\beta,1}\big)_{s}\leq
\frac{1}{\pi}\big\|\Psi_{\beta,n}(\cdot)\big\|_{s}\|
\varphi(\cdot)\|_{1}\leq\frac{1}{\pi}\big\|\Psi_{\beta,n}(\cdot)\big\|_{s}.
\end{equation}

Із співвідношень (\ref{yy}) i (\ref{q1}) за умов теореми 2 випливає
нерівність
\begin{equation}\label{t115}
\frac{1}{\pi}\|\Psi_{\beta,n}(t)\|_{s}\leq C_{a,b} \
(2s')^{\frac{1}{s}}\psi(n)\left(\eta(n)-n\right)^{\frac{1}{s}}, \ \
1<s\leq\infty, \ \frac{1}{s}+\frac{1}{s'}=1.
\end{equation}
Об'єднуючи (\ref{3t1}) i (\ref{t115}), одержуємо оцінку зверху для
величини ${\cal E}_{n}\big(L^{\psi}_{\beta,1})_{s}$ в співвідношенні
(\ref{theorem2}).

 Щоб одержати оцінку знизу величини
$E_{n}\big(L^{\psi}_{\beta,1}\big)_{s}$, $1<s\leq\infty$, розглянемо
функцію $f_{p}(t)$ вигляду (\ref{function}) при $p=1$, тобто функцію
вигляду
$$
f_{1}(t)=f_{1}(n, \psi, t)=
\frac{(a-1)(a-2)}{2\left(1+\pi^{2}\right)a(3a-4)}\times
$$
$$
\times\left(W_{[\eta(n)],[\eta(\eta(n))]}(\psi;0;t)-
W_{n,[\eta(n)]}(\psi;0;t)\right).
$$

В силу співвідношення (\ref{norm}) при $p=1$ для
$(\psi,\beta)$--похідної функції $f_{1}(\cdot)$ виконується
нерівність
$\Big\|\big(f_{1}(\cdot)\big)^{\psi}_{\beta}\Big\|_{1}\leq 1$, а
отже має місце включення
 ${f_{1}(\cdot)\in L^{\psi}_{\beta,1}}$.

Покажемо тепер, що при довільних
${\psi\in\mathfrak{M}^{+}_{\infty}}$,
${\lim\limits_{t\rightarrow\infty}(\eta(\psi,t)-t)=\infty}$,
${\eta(n)-n\geq a>2}, \ {\mu(n)\geq b>2}$, $\beta\in \mathbb{R}$
\begin{equation}\label{t15}
{E}_{n}(f_{1})_{s}\geq C_{a}\psi(n)(\eta(n)-n)^{\frac{1}{s'}}, \
1\leq s\leq\infty.
\end{equation}

З нерівностей (\ref{statement1}), (\ref{1t6}), (\ref{t9}) та
рівності (\ref{riv}), для будь-якого $t_{n-1}\in\mathcal{T}_{2n-1}$,
маємо
$$
\int\limits_{-\pi}^{\pi}\big(f_{1}(t)-t_{n-1}(t)\big)\left(W_{[\eta(n)],[\eta(\eta(n))]}(1;0;t)-
W_{n,[\eta(n)]}(1;0;t)\right)dt \leq
$$
$$
\leq\|f_{1}(t)-t_{n-1}(t)\|_{s}\left\|\left(W_{[\eta(n)],[\eta(\eta(n))]}(1;0;t)-
W_{n,[\eta(n)]}(1;0;t)\right)\right\|_{s'}=
$$
$$
=\|f_{1}(t)-t_{n-1}(t)\|_{s}\left\|\left(W_{[\eta(n)],[\eta(\eta(n))]}(\psi;0;t)-
W_{n,[\eta(n)]}(\psi;0;t)\right)^{\psi}_{0}\right\|_{s'}\leq
$$
\begin{equation}\label{3t8}
\leq
\frac{2\left(1+\pi^{2}\right)a(3a-4)}{(a-1)(a-2)}(\eta(n)-n)^{1-\frac{1}{s'}}
\|f_{1}(t)-t_{n-1}(t)\|_{s}.
\end{equation}

Зі співвідношень (\ref{t14}) i (\ref{3t8}) випливає нерівність,
справедлива для будь-якого $t_{n-1}\in\mathcal{T}_{2n-1}$
$$
\|f_{1}(t)-t_{n-1}(t)\|_{s }
\geq\frac{\pi}{96\left(1+\pi^{2}\right)^{2}}\frac{(a-1)^{2}(a-2)^{2}}{a^{3}(3a-4)}
\psi(n)(\eta(n)-n)^{\frac{1}{s'}}=
$$
\begin{equation}\label{3t9}
=C_{a,b}\psi(n)(\eta(n)-n)^{\frac{1}{s'}}.
\end{equation}
Тим самим нерівність (\ref{t15}) доведено, а разом з нею і теорему
2.

\textbf{Наслідок 2.} \emph{ Нехай $\psi_{r}(t)=\exp\left(-\alpha
t^{r}\right), \  {r\in(0,1)}, \ {\alpha>0}$, ${1< s\leq \infty}$,
$\frac{1}{s}+\frac{1}{s'}=1$, $\beta\in \mathbb{R}$. Тоді  для всіх
номерів $n$, таких, що $$n\geq
1+\max\left\{\left(\frac{2r\alpha}{\ln2}\right)^{\frac{1}{1-r}},
2\left(\frac{\ln2}{\alpha\left(3^{r}-2^{r}\right)}\right)^{\frac{1}{r}}\right\},$$
 справедливі оцінки}
$$
C_{a} \exp\left(-\alpha
n^{r}\right)n^{\frac{1}{s'}}\left(\left(\frac{\ln2}{\alpha
n^{r}}+1\right)^{\frac{1}{r}}-1\right)^{\frac{1}{s'}} \leq \ \ \ \ \
\ \ \ \ \ \ \ \ \ \ \ \ \ \ \ \ \ \ \ \ \ \ \ \ \ \ \
$$
$$
\leq{ E}_{n}\Big(L^{\psi_{r}}_{\beta,1}\Big)_{s}\leq{\cal
E}_{n}\big(L^{\psi_{r}}_{\beta,1})_{s}\leq
$$
$$
\ \ \ \ \ \ \ \ \ \ \ \ \ \ \ \ \ \ \ \ \ \ \ \ \ \ \ \ \ \  \leq
C_{a,b} \ (2s')^{\frac{1}{s}} \exp\left(-\alpha
n^{r}\right)n^{\frac{1}{s'}}\left(\left(\frac{\ln2}{\alpha
n^{r}}+1\right)^{\frac{1}{r}}-1\right)^{\frac{1}{s'}},
$$
\emph {де величини $C_{a}$ i $C_{a,b}$ означаються формулами
(\ref{Ca}) i (\ref{Cab}) при ${a=a(\alpha,r)}$, $b=b(\alpha,r)$, що
в свою чергу задані за допомогою рівностей (\ref{aa}) i (\ref{bb})
відповідно.}

Також для величин  $E_{n}\big(L^{\psi_{r}}_{\beta,1}\big)_{s}$,
$r\in(0,1], \ {\alpha>0}$, ${1< s\leq \infty}$,
${\frac{1}{s}+\frac{1}{s'}=1}$, можна записати аналогічне до
(\ref{consequence}) співвідношення
$$
\mathcal{E}_{n}\Big(L^{\psi_{r}}_{\beta,1}\Big)_{s}\asymp{
E}_{n}\Big(L^{\psi_{r}}_{\beta,1}\Big)_{s}\asymp\exp\left(-\alpha
n^{r}\right)n^{\frac{1-r}{s'}}, \ \ n\in\mathbb{N}.
$$

E-mail: \href{mailto:serdyuk@imath.kiev.ua}{serdyuk@imath.kiev.ua},
\href{mailto:tania_stepaniuk@ukr.net}{tania$_{-}$stepaniuk@ukr.net}

\end{document}